\documentclass[10pt,psamsfonts]{amsart}

\makeatletter

\@addtoreset{equation}{section}
\makeatother

\usepackage{epsfig}
\usepackage{graphicx}
\usepackage{amsmath}
\usepackage{amssymb}
\usepackage[latin1]{inputenc}
\usepackage{float}
\usepackage{url}
\RequirePackage{amssymb}
\RequirePackage{latexsym}
\RequirePackage{subeqnarray}

\begin{document}
\thispagestyle{empty} \setcounter{page}{1}

\title[Investigation of an endogenous reproductive rhythm]{Endogenous circannual rhythm in LH secretion: insight from signal analysis coupled with mathematical modelling}

\author[A. Vidal, C. Médigue, B. Malpaux, F. Clément]{Alexandre Vidal$^{1}$ , Claire M\'edigue$^{1}$ , Beno\^it Malpaux$^{2}$  and Fr\'ed\'erique Cl\'ement$^{1}$}

\maketitle

\begin{center}
$^{1}$ INRIA Paris-Rocquencourt Research Centre, Project-Team SISYPHE.\newline
\newline
$^{2}$ INRA-CNRS-Université de Tours-Haras Nationaux,

UMR Physiologie de la Reproduction et des Comportements.
\end{center}

\footnotetext[1]{Domaine de Voluceau, Rocquencourt BP 105, 78153 Le Chesnay cedex, FRANCE.}
\footnotetext[2]{Centre INRA de Tours, 37380 Nouzilly, FRANCE.}

\begin{abstract}

In sheep as in many vertebrates, the seasonal pattern of reproduction is timed by the annual photoperiodic cycle, characterized by seasonal changes in the daylength. The photoperiodic information is translated into a circadian profile of melatonin secretion. After multiple neuronal relays (within the hypothalamus), melatonin impacts GnRH (gonadotrophin releasing hormone) secretion that in turn controls ovarian cyclicity. The pattern of GnRH secretion is mirrored into that of LH (luteinizing hormone) secretion, whose plasmatic level can be easily measured. We addressed the question of whether there exists an endogenous circannual rhythm in a tropical sheep (Black-belly) population that exhibits clear seasonal ovarian activity when ewes are subjected to temperate latitudes. We based our analysis on LH time series collected in the course of 3 years from ewes subjected to a constant photoperiodic regime. Due to intra- and inter- animal variability and unequal sampling times, the existence of an endogenous rhythm is not straightforward. We have used time-frequency signal processing methods, and especially the Smooth Pseudo-Wigner-Ville Distribution, to extract possible hidden rhythms from the data. To further investigate the LF (low frequency) and HF (high frequency) components of the signals, we have designed a simple mathematical model of LH plasmatic level accounting for the effect of experimental sampling times. The model enables us (i) to confirm the existence of an endogenous circannual rhythm as detected by the LF signal component, (ii) to investigate the action mechanism of photoperiod on the pulsatile pattern of LH secretion (control of the interpulse interval) and (iii) to conclude that the HF component is mainly due to the experimental sampling protocol.

\vspace{3mm}


\end{abstract}

\section{Introduction}
In temperate regions, sheep exhibit annual cycles of breeding activity (see \cite{Karsch89}, \cite{Chemineau07}), timed by the photoperiodic cycle. The photoperiodic information is translated into circadian profiles of melatonin secretion. After multiple neuronal relays (within the hypothalamus), melatonin impacts GnRH (gonadotrophin releasing hormone) secretion that in turn controls ovarian cyclicity. The pattern of GnRH secretion is mirrored into that of LH (luteinizing hormone) secretion, whose plasmatic level can be easily measured. It has been generally assumed that the photoperiod entrains a circannual endogenous rhythm of sexual activity to a period of one year. This assumption has not been yet verified in sheep native from tropical regions. Although previous studies have evidenced seasonal changes in LH secretion in these ewes, it remains unknown whether such changes reflect an endogenous rhythm (see \cite{Chemineau04}).

This study aims at gathering evidence for the existence of an endogenous circannual reproductive rhythm in Black Belly ewes, from the study of long-term variations in LH plasmatic levels. To get rid of the entrainment by the photoperiodic cycle, the ewes were exposed to a fixed day length corresponding to permanent short days.

\begin{figure}[htbp]
\centering
\includegraphics[height=18cm]{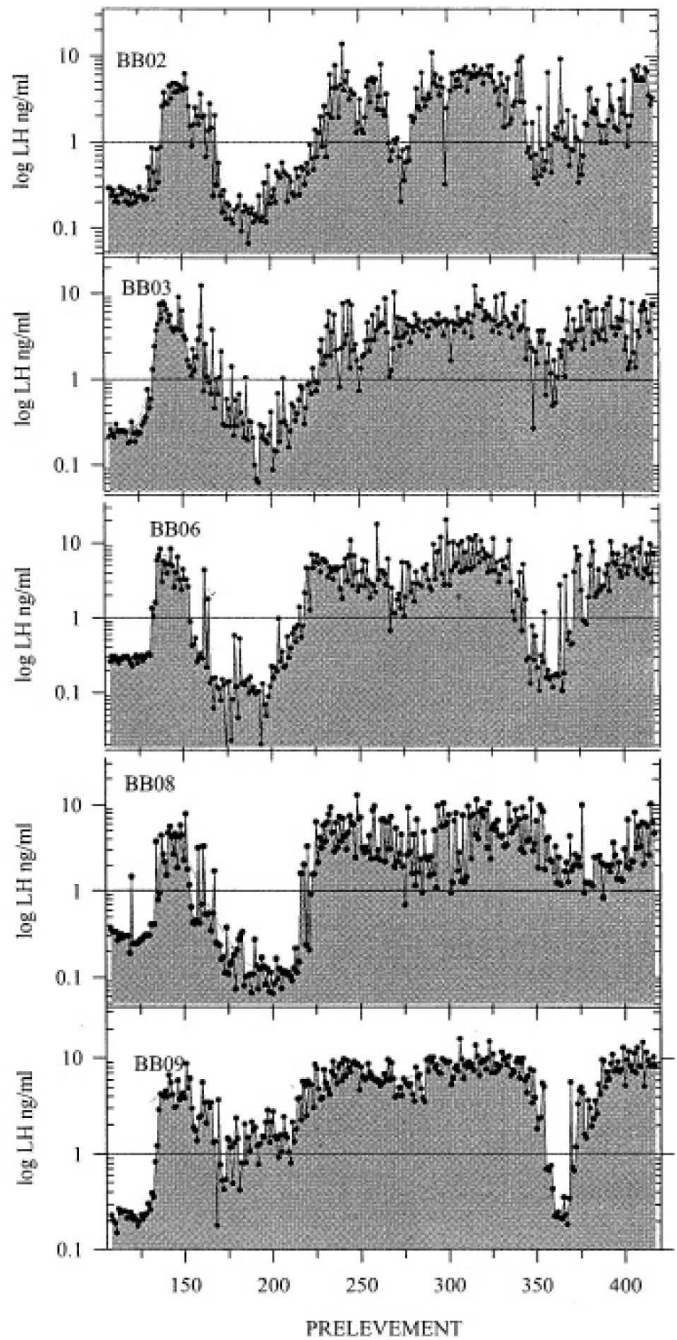}
\caption{ {\em  Changes in the amplitude of LH plasmatic levels for the five ewes, exposed to a short day regimen for three years.}
All ewes  exhibit a low signal-to-noise ratio, due to a strong high frequency (HF) activity, in the range of the sampling frequency of the measure, which complicates the research for an hidden endogenous circannual rhythm.}
\label{fig1}
\end{figure}

The raw LH time series are characterized by a low signal-to-noise ratio (see figure \ref{fig1}). Since extracting a circannual trend from the raw data amounts to detect time-varying hidden rhythms, we chose to use a time-frequency method, the Smooth Pseudo Wigner-Ville Distribution (SPWVD). This method provides meaningful time-varying spectral parameters and can discriminate efficiently the low frequency component (LF) against the high frequency component (HF) in the LH signals. The SPWVD has enabled us to reveal a clear endogenous rhythm from the LF component [and to compare its properties to those of a purely periodic annual rhythm]. It has also revealed that the HF component represents a well located and quantifiable activity different from a white noise.

To further investigate the LF and HF components of the signals, we have designed a simple mathematical model of LH plasmatic level, accounting for the main mechanisms underlying LH secretion and the effect of experimental sampling times. The comparative SPWVD analysis of the model outputs with the experimental data has enabled us (i) to confirm the existence of an endogenous circannual rhythm as detected by the LF signal component, (ii) to investigate the action mechanism of photoperiod on the pulsatile pattern of LH secretion (control of the interpulse interval) and (iii) to conclude that the HF component is mainly due to the experimental sampling protocol.

\section{Material and Methods}

\subsection{Animals and experimental design}

The study was carried out with five Black-Belly ewes during three years, from January 2002 to January 2005. In January 2002, the ewes were ovariectomized and received a constant-release estradiol implant (OVX+E$_{2}$). They were exposed to an alternating 8 hour light/16 hour dark regime (8L:16D) mimicking a short day length. The subcutaneous Silastic implant contains a 30 mm column of \oe{}stradiol-17$\beta$ to maintain a serum concentration of 3-5 pg ml$^{-1}$, which is intermediate between the concentrations observed in the luteal and follicular phase of the ovarian cycle. The OVX+E$_{2}$ ewes were exposed for three months to an inhibitory long photoperiod (16L:8D), prior to their entering the short day regime (LD-SD transition). They were kept in the light-controlled room of an experimental breeding station near Tours (France, 47$^{0}$ N). The artificial photoperiodic regime was tuned by means of electric timers, whereas temperature was not controlled. The reproductive endocrine state was monitored by regular measurements of LH plasmatic level. Blood samples were collected from the jugular vein, twice a week, in the beginning of the morning and stored at $-20^{0}$ C. The sensitivity of the plasma LH assay was 0.2 ng ml$^{-1}$.

\subsection{Data Analysis : Signal Processing}

To follow instantaneous changes in spectral activities, we selected as time-frequency method the Smooth Pseudo Wigner-Ville Distribution (SPWVD), which is based on the Wigner-Ville distribution (see \cite{Cohen89}, \cite{Pola96}). The Wigner-Ville distribution is always real-valued, with a very good time and frequency resolution preserving time and frequency shifts. This distribution satisfies the marginal properties: the energy spectral density and the instantaneous power can be obtained as marginal distributions of the Wigner-Ville distribution. As a main drawback, interference between components of different frequency may occur, leading to oscillatory cross-terms (see \cite{Martin85}). To overcome this problem, the SPWVD is computed as a time average of the Wigner-Ville distribution: 
\begin{equation}\label{eqn:SPWV_def}
\begin{array}{ll}
SPW(t,\nu)=\int_{-\infty}^{+\infty}\int_{-\infty}^{+\infty}g_M(s-t)h^2_N(\tau)x(s+\tau/2)\cdot\\
x^{*}(s-\tau/2)\,e^{-2i \pi \nu \tau}\ d\tau \ ds
\end{array}
\end{equation}
where N and M are respectively the size of the $h_N$ and $g_M$ windows, the $h_N$ window weights the input time-series $X(t)$ and the $g_M$ window averages the instantaneous spectrum. In the SPWVD, a trade-off is required between the time-frequency resolution and the removal of the cross-terms: smoothing in time by $g_M$ and/or in frequency by $h_N$ reduces the time resolution and/or the frequency resolution (see \cite{Pola96}).

Even if the marginal properties are no longer strictly valid and the time resolution reduced (depending on the duration of $g_M$), the SPWVD can be used to estimate both the instantaneous power ($IPow$) and frequency ($IF$) of the signal from the first two order moments of the distribution:
\begin{equation}
\label{eqn:SPWVD_mom}
\begin{array}{lll}
IPow(t) &=& \int\limits_{-\infty }^{+\infty }SPW(t,\nu)\ d\nu \\
IF(t)   &=& \frac{ \int\limits_{-\infty}^{+\infty} \nu\ SPW(t,\nu)\ d\nu}{ \int\limits_{-\infty}^{+\infty} SPW(t,\nu)\ d\nu}
\end{array}
\end{equation}
The instantaneous power ($IPow$) can be converted to the instantaneous amplitude as follows:
\begin{equation}\label{eqn:SPWVD_IAmp}
IAmp(t)=\alpha\ \sqrt{IPow(t)}
\end{equation}
where $\alpha$ is a constant depending on the filter properties.

The time-frequency analysis was performed within the Scilab-Scicos academic environment\footnote{\url{http://www.scilab.org/} \url{http://www.scicos.org/}.}.
Due to interference between the spectra, the use of SPWVD needs some preprocessing of the LH time series. The data were resampled by interpolation with a third-order spline function to get equally spaced time series with a daily sampling time. Two mono frequency components were extracted after filtering. The LF component corresponds to long-term fluctuations in LH level (alternation between ovarian cyclic activity and inactivity), while the HF component covers the range of the experimental sampling frequency. Accordingly, the LF activity was upper-bounded at 5 months with a low-pass finite impulse response filter (FIR). The 5-month bound was chosen to encompass relatively short term changes in LH levels. The HF activity was bounded between 1.5 and 3 weeks (around the experimental sampling frequency) with a band-pass FIR.

\begin{figure}[htbp]
\centering
\includegraphics[width=8cm]{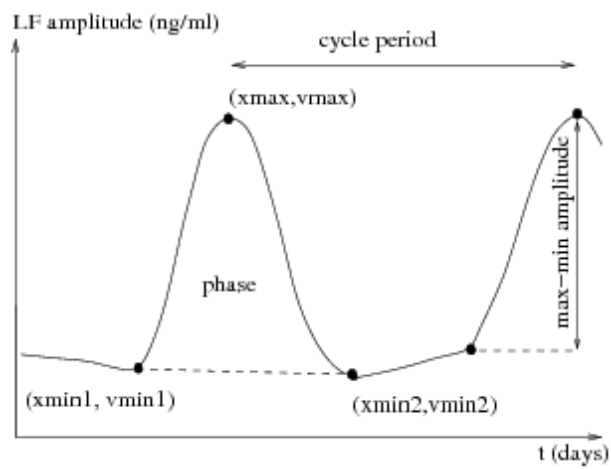}
\caption{ {\em LF parameters extracted from SPWVD  to characterize the phases of  LH secretion.}
xmax (day), vmax (ng ml$^{-1}$): occurrence and value of the maximal LF amplitude of a phase; xmin (day), vmin (ng ml$^{-1}$): occurrence and value of the minimal LF amplitude of a phase; max-min amplitude(ng ml$^{-1}$) : difference between the maximal and minimal LF values of a phase; mean amplitude (ng ml$^{-1}$ day$^{-1}$ ): sum of instantaneous amplitudes  normalized by the number of days between xmin1 and xmax1;  phase duration (day) : duration between the midpoint joining  xmin1 and  xmax, and the midpoint joining  xmax and  xmin2 of a phase;  cycle period (day): duration between two xmax.}
\label{fig2}
\end{figure}

After SPWVD processing, several parameters were extracted from the resulting time series for LF amplitude (as illustrated in figure \ref{fig2}), to describe phases and cycle periods of LH secretion. An activity phase encompasses the period where plasmatic LH levels are high, resulting from the high frequency LH spikes in the sexual season. An activity phase is defined as a time interval $[t_1,t_2]$ such that:
\begin{enumerate}
\item the LH level is over 0.2 ng ml$^{-1}$ during this whole interval,
\item the signal derivative at time $t_1$ becomes greater than a threshold (7.10$^{-3}$ ng ml$^{-1}$ day$^{-1}$),
\item the upper bound $t_2$ is the lowest time greater than $t_1$ at which the signal derivative is negative, increasing and becomes greater -7.10$^{-3}$ ng ml$^{-1}$ day$^{-1}$.
\end{enumerate}
A cycle period corresponds to the sequence of one activity phase and the subsequent inactivity phase.
\begin{itemize}
 \item xmax (day) and vmax (ng ml$^{-1}$) stand for the occurrence and value of the maximal amplitude of a phase; hence the time separating two subsequent values of xmax is analogue to a cycle period.
 \item xmin (day) and vmin (ng ml$^{-1}$) stand for the occurrence and value of the minimal amplitude of a phase; 
 \item max-min amplitude (ng ml$^{-1}$) stands for the difference between the maximal and minimal values of a phase; 
 \item mean amplitude (ng ml$^{-1}$ day$^{-1}$) stands for the sum of instantaneous amplitudes normalized by the number of days between xmin1 and xmax1;
 \item phase duration (day) begins at the midpoint joining xmin1 and xmax, and ends up at the midpoint joining xmax and xmin2 of a phase.
\end{itemize}
These parameters are useful in assessing the existence of a circannual trend in the LH series and comparing the LH series between the 5 ewes of the study. In particular, the synchronisation is estimated by the standard error to mean (SEM) of the time occurrences of the minimal and maximal amplitude values in each phase. Average parameter values can also be used to compare the Black Belly strain to other strains originating from temperate regions such as the Suffolk strain. Statistical comparisons have been performed as a one way repeated measure analysis of variance (anova). All results are expressed as mean and SEM.

\subsection{Mathematical modelling}

To assess the results of the time-frequency analysis and improve our interpretation of the frequency-filtered components, we have built a simple mathematical model of the plasmatic LH level. It is based on a representative function of the pulsatile LH secretion by the pituitary gland coupled with a term accounting for dissemination into the blood. The LH pulse frequency is subject to the influence of the photoperiodic regime. More precisely, basing on available knowledge, we have introduced the influence of the light/dark regime on the LH secretion rate at the pituitary level. To compare the photoperiod-entrained circannual rhythm with the free-running rhythm, we have tested the effect of both sinusoidal and damped-sinusoidal terms. Synthetic (model output) LH series comparable to the experimental series have been obtained after mimicking the process of experimental sampling and have been subjected to the same SPWVD analysis as the experimental data.

\subsubsection{LH secretion by the pituitary gland}

\begin{figure}[htbp]
\centering
\includegraphics[width=12cm]{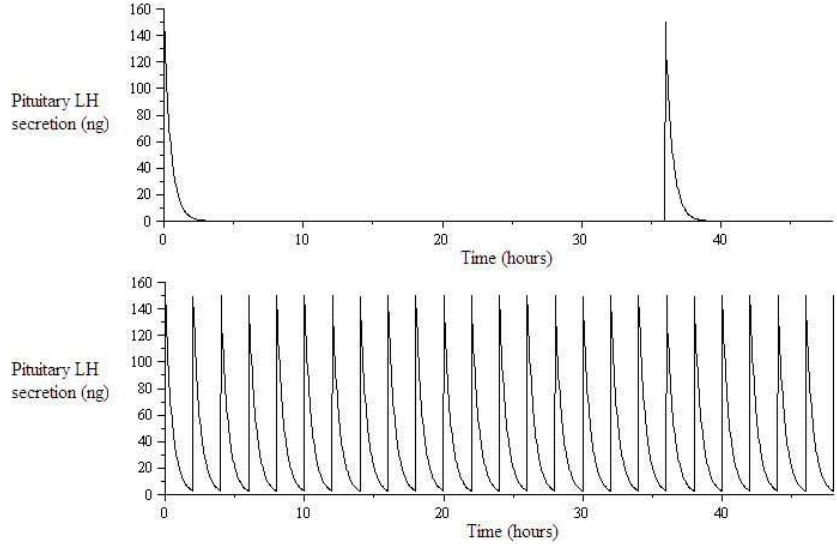}
\caption{ {\em Graphs $LH(t)$ when $P_{spike}(t)\in \{2,36\}$.}
The top (resp. bottom) panel represents $LH(t)$ on a 48 hours time interval when $P_{spike}(t)$ is constant and equals 36 (resp. 2).}
\label{fig3}
\end{figure}

The pituitary gland releases LH as successive spikes. We chose to represent each spike of LH release by an instantaneous impulse followed by a fast exponential decrease. The interspike interval is controlled by a time varying function (accounting for the photoperiod action) $P_{spike}(t)$. Hence, the instantaneous release of LH by the pituitary gland is given by:
\begin{equation}
LH(t)=A_{spike}\exp \left[ -k_{hl}\left( t-\left\lfloor \frac{t}{P_{spike}(t)}\right\rfloor P_{spike}(t)\right) \right]
\label{LHspike}
\end{equation}
where $\lfloor x \rfloor$ refers to the greatest integer smaller than $x$ (integer part). Figure \ref{fig3} displays the graph of the $LH(t)$ function when $P_{spike}(t)$ is constant and equals 36 (top panel) and 2 (bottom panel). We will give the precise expressions for $P_{spike}(t)$ in the next subsection. The $A_{spike}$ parameter represents the spike amplitude. The $k_{hl}$ stiffness coefficient is directly linked to the spike half-life, $\tau _{hl}$, by the following relation:
\[
k_{hl}=\frac{\ln 2}{\tau _{hl}}
\]

The parameter values are chosen so as to fit at best the shape of the LH pulse in the ewe. We fist base on the few experiments that have estimated the LH release by the pituitary gland in the ewe from regular sampling in portal hypophyseal blood (see \cite{Pincus98}, \cite{Keenan04}). We also use the results from application of deconvolution methods to time series of LH sampled in jugular blood that give indications about the amplitude and the half-life of the LH spike (see \cite{Veldhuis87}, \cite{Veldhuis90}, \cite{DeNicolao99}). Accordingly, we chose the following parameter values for the model outputs : $A_{spike}=150$ ng and $k_{hl}=2$, to obtain a spike half-life $\tau _{hl}$ around $20$ minutes.

\subsubsection{Control of the interspike interval}

The $P_{spike}$ function represents, in a compact way, the effect of light on LH secretion at the pituitary level and reproduces the alternation between high LH release frequency (i.e. short interspike interval) and low LH release frequency (i.e. long interspike interval). We will use two different expressions of $P_{spike}$ to mimic either the natural temperate photoperiod or the constant SD light regime used in the experimental protocol.

Information on the daylight duration is encoded as melatonin secretion by the pineal gland and conveyed to GnRH neurons through complex neuronal circuits. GnRH secretion controls in turn the frequency of LH secretion. In ewes, LH pulse frequency increases as the daylight duration gets shorter, which sets sexual activity on. In contrast, LH pulse frequency decreases as the daylight duration gets longer, leading to sexual inactivity. Approximating the annual variations in daylight by a sinusoidal function, we can introduce the impact of the light regime in temperate regions as sustained oscillations of the interspike interval:
\begin{equation}
P_{spike}(t)=\frac{\Delta P}{2}\left( \sin \left( \frac{2\pi t}{P_{photo}}\right) +1\right) +P_{spike}^{\min } 
\label{PSin}
\end{equation}
\begin{equation*}
\Delta P=P_{spike}^{\max }-P_{spike}^{\min } 
\end{equation*}
$P_{spike}^{\max }$ corresponds to the maximal interspike interval during the inactivity phase in the ewe ($36$ hours), $P_{spike}^{\min }$ corresponds to the minimal interspike interval in the activity phase ($1.5$ hours) and the photoperiod $P_{photo}$ is one year ($8760$ hours).

In the experimental constant SD light regime, we may expect that the oscillations in $P_{spike}$ progressively weaken since there is no more entrainment.
This can be tested using a damped sinusoid for $P_{spike}$, for which the annual increase in the interspike interval is less and less pronounced with time :
\begin{equation}
P_{spike}(t)=\left\{ 
\begin{array}{lcc}
\frac{\Delta P}{2}\left( \sin \left( \omega t\right) +1\right)+P_{spike}^{\min } & \text{if} & t<P_{photo} \\ 
\frac{\Delta P}{2}\left( \sin \left( \sqrt{1-\zeta ^{2}} \omega t\right) +1 \right) K e^{-\zeta \omega (t-P_{photo})}+P_{spike}^{\min } & \text{if} & t\geq P_{photo}
\end{array}
\right.
\label{PSinAm}
\end{equation}
where :
\[
\omega=\frac{2 \pi}{P_{photo}}
\]
For $t$ in $[0, P_{photo}]$ the expression of $P_{spike}$ is the same as before. For $t>P_{photo}$, the sinusoid is damped with coefficient $\zeta $ ($\zeta =0.14$).
For sake of continuity of $P_{spike}$, $K$ is given by:
\[
K=\frac {1}{\left( \sin \left( 2\,\sqrt {1-{\zeta }^{2}}\pi \right) +1 \right) }
\]
Figure 4 illustrates the changes in $P_{spike}(t)$ with the sinusoidal $P_{spike}$ function defined by (\ref{PSin}) (top panel) and the damped sinusoidal $P_{spike}$ function defined by (\ref{PSinAm}) (bottom panel) with:
\begin{equation}
\begin{array}{c}
P_{spike}^{\max }=36\text{ hours}\quad ;\quad P_{spike}^{\min }=1.5\text{ hours}\quad ;\quad \Delta P=34.5\text{ hours} \\
P_{photo}=8760\text{ hours}
\end{array}
\label{ParamSin}
\end{equation}

\begin{figure}[htbp]
\centering
\includegraphics[width=9cm]{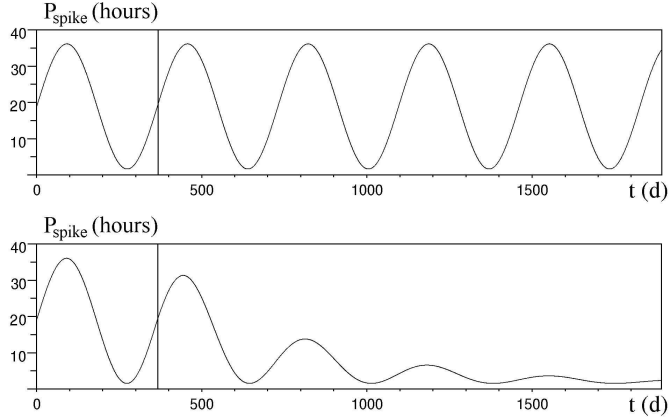}
\caption{ {\em Graphs of the two functions used successively for $P_{spike}$.}
The top panel represents the sinusoid function defined by (\ref{PSin}); the second one is defined by (\ref{PSinAm}) as the same sinusoid over $[0, 1\ year]$, and, over $[1\ year, 5\ years]$, a damped sinusoid with same fundamental period and damping coefficient $\zeta=0.14$. The vertical segment at 1 year represents the beginning of the sampling process on the fine numerical integration of equation (\ref{LHplasma}).}
\label{fig4}
\end{figure}

\subsubsection{Plasmatic LH level and synthetic sampling process}

From (\ref{LHspike}), LH blood level changes according to :
\begin{equation}
\frac{dLH_{p}}{dt}(t)=LH(t)-\alpha LH_{p}(t)
\label{LHplasma}
\end{equation}
where $\alpha$ represents the LH clearance rate from the blood. To be consistent with the one hour half-life of LH pulses in the jugular blood, we have fixed $\alpha =6$.

Then, we have extracted (unequally spaced) time-series from the fine step numerical integration of equation (\ref{LHplasma}), with a sampling frequency of twice a week (comparable to the experimental protocol). The sampling times (expressed in hours) are defined by:
\[
t_{i}=84t+a(i)
\]
The random numbers $a(i)\in [-0.5,0.5]$ are generated through a Mersenne Twister algorithm (see \cite{Matsumoto98}) to take into account the inherent variability of the sampling times.

Once we have simulated the model outputs, we dispose of synthetic time series that can be subject to the same SPWVD analysis as the experimental data. We will take advantage of our full control of the frequency components to get insight on the interpretation of the LF and HF components of the experimental time series. In particular, we will be able to study the HF component and separate the variability due to the experimental protocol.

\section{Results}
\subsection{Signal analysis}
\subsubsection{Low frequency band - LH secretion phases}
\begin{figure}[htbp]
\centering
\includegraphics[width=9.5cm]{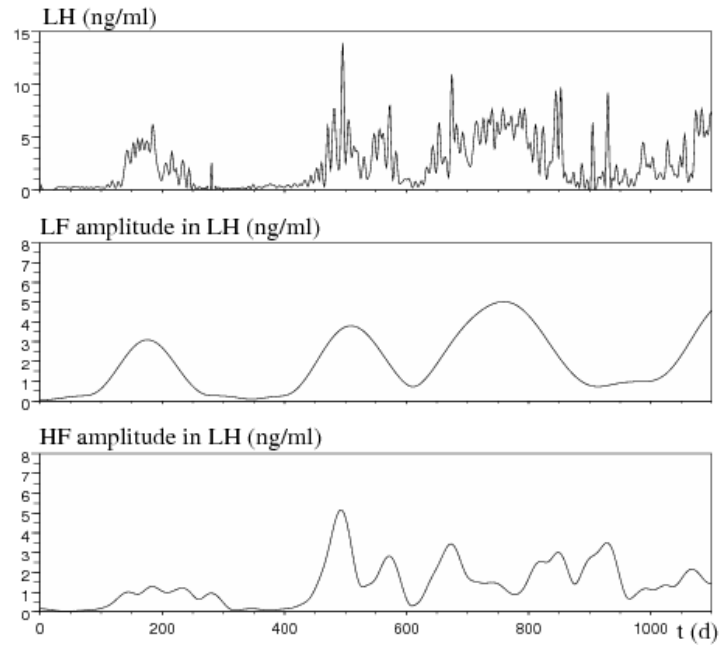}
\caption{ {\em Time evolution of the LH plasmatic secretion and its low and high frequency components during the three years of the experiment, in a ewe.}
The top panel represents the raw LH time series, the middle panel and the bottom panel successively represents the LF and HF amplitude time series, extracted from the raw series by the SPWVD analysis.}
\label{fig5}
\end{figure}

Figure \ref{fig5} shows the time changes in the LH plasmatic secretion along the three years of the experiment, in one ewe. The top panel, represents the raw LH time series. The middle panel represents the LF amplitude time series, after being extracted from the raw series by the SPWVD analysis. The LF amplitude captures the basal LH plasmatic secretion correctly and its profile does appear to be compatible with an endogenous circannual rhythm.

Four phases have been evidenced in each ewe during the three years, the last phase being truncated. The secretion parameters are summarized in table \ref{tab1}, as mean values (averaged over the five ewes) for each phase of LH plasmatic secretion. These results show a significant change in the LH secretion in the course of the experiment: the minimal and maximal values, as well as the duration of the LH phases progressively increase from phase one to three, so that the phase mean amplitude gets greater.

\begin{figure}[htbp]
\centering
\includegraphics[width=10.5cm]{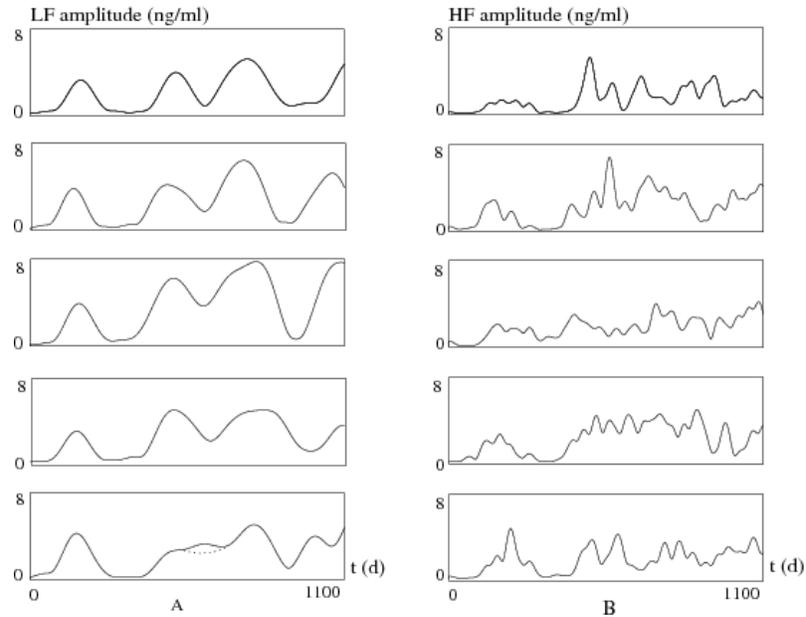}
\caption{ {\em  Time evolution of the amplitude of the low and high frequency components of the  LH plasmatic time series for the five ewes.}
Panel A shows the LF amplitude, reflecting the basal LH plasmatic secretion. Its time evolution, which reveals four phases, is compatible with an endogenous circannual rhythm. In each ewe there is a progressive increase in the basal level of LH plasmatic secretion. Panel B represents the HF amplitude, which increases during the experiment, especially from the second year on.}
\label{fig6}
\end{figure}

\begin{figure}[H]
\centering
\includegraphics[width=8.5cm]{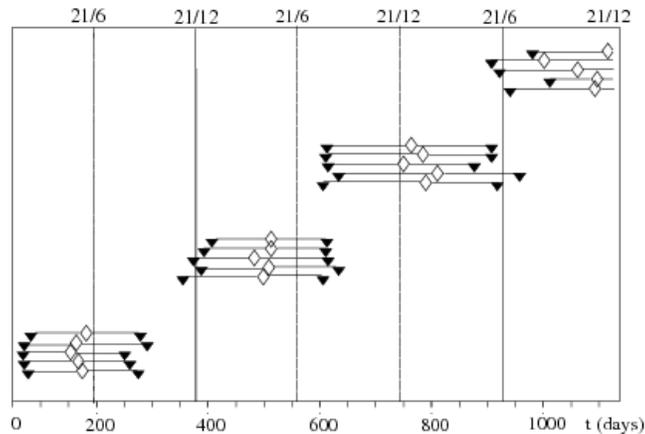}
\caption{{\em Timing of LH phases in the five Black-Belly ewes during three years.}
Black triangles mark the beginning and end of a LH phase of secretion; white diamonds mark for the maximal amplitude of the phase. The degree of synchronization decreases, while the duration of the phases increases from the first to the last phase. The latest two cycle periods (duration between two maximal amplitude values of a phase) are shorter.}
\label{fig7}
\end{figure}

\begin{table}[hb]
\caption{Phases of LH secretion: LF parameters values}
\begin{center}
\begin{tabular}{|l|l|l|l|l|}
\hline
LF parameters& 
Phase 1 & 
Phase 2 & 
Phase 3 & 
p \\
\hline
xmin, xmax synchro (day) & 
2 ; 4& 
8 ; 5& 
4 ; 10& 
 \\
\hline
min. amp (ng/ml)& 
0.37 $\pm $ 0.06& 
0.42 $\pm $ 0.08& 
2.25$\pm $ 0.52& 
0.003 \\
\hline
max. amp (ng/ml)& 
3.55$\pm $ 0.24& 
4.30$\pm $ 0.59& 
5.77$\pm $ 0.5& 
0.01 \\
\hline
mean amp (ng/ml/day)& 
1.96$\pm $ 0.13& 
2.88$\pm $ 0.36& 
3.95$\pm $ 0.43& 
0.001 \\
\hline
max-min amp(ng/ml)& 
3.18$\pm $ 0.19& 
3.88$\pm $ 0.51& 
3.52$\pm $ 0.53& 
0.3 \\
\hline
phase duration (day)& 
98$\pm $3& 
116$\pm $ 4& 
149$\pm $ 5& 
0.001 \\
\hline
\end{tabular}
\label{tab1}
\end{center}
\end{table}

Figure \ref{fig6} (panel A) and \ref{fig7} illustrate the individual change in the LF amplitude for each of the 5 ewes. Figure \ref{fig6}, panel A, illustrates the progressive increase in the basal level of LH. Figure \ref{fig7} displays the characteristic times (beginning, end and maximal value) of the LH secretion phases: the synchronisation between the ewes decreased whereas the phase duration increased from the first to the last phase. The cycle periods were successively 334$\pm $ 4, 276$\pm $9, 294 $\pm $22, with p $\leq $ 0.07. It is worth noticing that the decrease in the second cycle period corresponds to a resetting of the third LH phase secretion according to the temperature cycle, which persists during the fourth phase. Hence, the times of maximal LF amplitude coincide with the time of minimal temperature values.

\subsubsection{High frequency band }
Figure \ref{fig5} shows the change in LH plasmatic secretion during the three years of the experiment, in one ewe. The top panel represents the raw LH time series, the middle and bottom panel respectively represents the LF and HF amplitude time series, after being extracted from the raw series by the SPWVD analysis. The level of HF amplitude increases during the experiment, especially from the second year on. The HF amplitude is not negligible at all compared to the LF one, which is consistent with the difficulty encountered to elucidate a LF rhythm from raw series, as shown on figure \ref{fig1}. It is worth noticing that the HF and LF peaks do not occur simultaneously. It is worth noticing that the HF and LF peaks do not occur simultaneously. This tendency is shared between the 5 ewes (figure \ref{fig6}, panel B).

\subsection{Mathematical modelling}

\subsubsection{A sinusoidal evolution of the interspike interval models the influence of natural temperate photoperiod}

\begin{figure}[htbp]
\centering
\includegraphics[width=11cm]{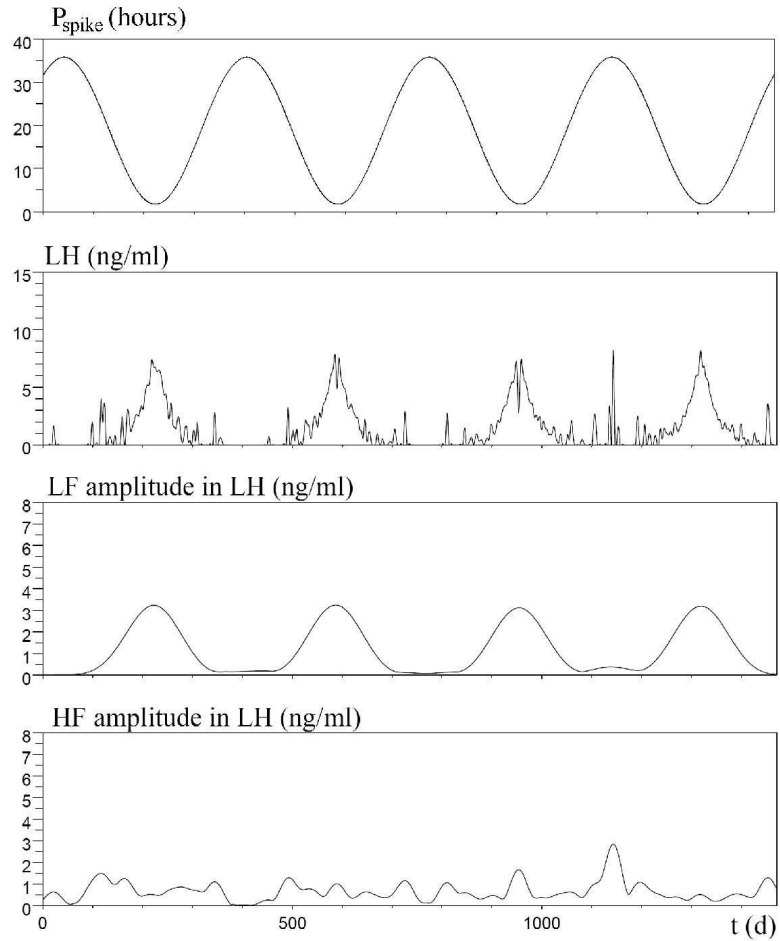}
\caption{ {\em Model output for a sinusoidal interspike interval function and its low and high frequency components during four years.}
The first panel represents the sinusoid chosen for $P_{spike}$ to model the answer to natural temperate photoperiod. It is drawn from the beginning of sampling process.
The second panel represents a raw LH time series obtained from the sampling process (twice per week) upon a model output. The third and fourth panels successively represent the LF and HF amplitude time series, extracted from the raw series by the SPWVD analysis.}
\label{fig8}
\end{figure}

The second panel of figure \ref{fig8} displays a typical output of the model with the sinusoidal function $P_{spike}$ defined in (\ref{PSin}). When $P_{spike}$ takes its highest values, LH pulses are scarce, so that the sampling times occur most of the time when the plasmatic LH level is negligible (under the detection threshold 0.2 ng ml$^{-1}$). However, it may happen that a sampling time occurs shortly after a spike, leading to isolated great sampled values such as that visible around 1150 days. When the values of $P_{spike}$) are at their lowest, the LH pulses are much more frequent, so that the basal level remains at a relatively elevated level between two pulses, while the maximal values achieved after a spike increase. Moreover, the basal level of $LH(t)$ increases leading necessarily to great sample values.

The LF component, extracted from SPWVD analysis, distinguishes between the activity and inactivity phases (see third panel of figure \ref{fig8}).
The LF maxima all have the same amplitude, which is a characteristic feature of the response to the natural temperate photoperiod.

The HF amplitude is not negligible but remains very weak compared to that of LF (see fourth panel of figure \ref{fig8}). Furthermore, there is no systematic coincidence between the changes in HF activity and LF activity. It is also worth noticing that the HF component does not impact on the LF component dynamics in a significant way. For instance, despite the sequence of great spikes around 1150 days (during an inactivity phase), the LF activity pattern just undergoes a small variation (less than 8\% of the maximal amplitude).

Finally, the 4 to 6 months activity phases are similar to the duration of sexual activity phases in seasonal mammals under natural photoperiod.

\subsubsection{A damped sinusoidal evolution of the interspike interval models the LH sampled level under LD-SD light control}

\begin{figure}[H]
\centering
\includegraphics[width=10.45cm]{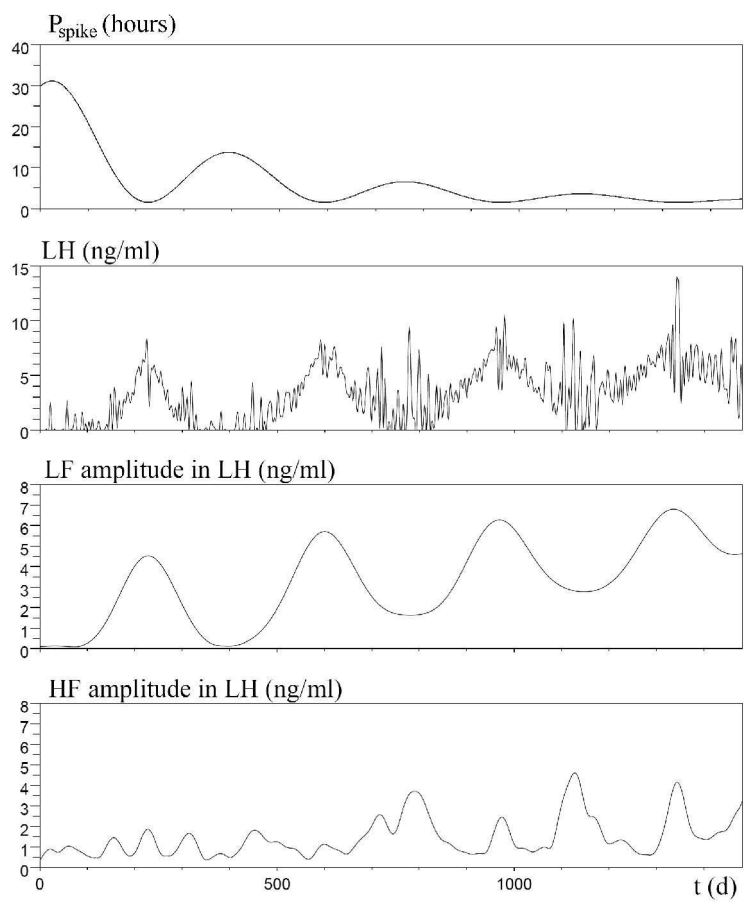}
\caption{ {\em Model output for a damped sinusoidal interspike interval function and its low and high frequency components during four years.}
The first panel represents the damped sinusoid chosen for $P_{spike}$ to model the answer to SD light control. It is drawn from the beginning of sampling process.
The second panel represents a raw LH time series obtained from the sampling process (twice per week) upon a model output. The third and fourth panels successively represent the LF and HF amplitude time series, extracted from the raw series by the SPWVD analysis.}
\label{fig9}
\end{figure}

The second panel of figure \ref{fig9} displays a typical output of the model with the damped sinusoid $P_{spike}$ function defined in (\ref{PSinAm}). We can easily identify a first activity period similar to that obtained with a sinusoid function, which consists of a few spikes during the inactivity phase followed by an increase in LH level and a return to a zero basal level. After the first year, the transitions between the activity and inactivity phases are less obvious, since the signal becomes noisier.

The LF component extracted from SPWVD analysis (see third panel of figure \ref{fig9}) displays 3 other phases, one each year. The maximal LF amplitude (occurring when $P_{spike}(t)$ is minimal) increases from one cycle to the next, while the basal level (corresponding to the local maxima of $P_{spike}(t)$) increases from one year to the next.

The HF activity (displayed in the 4th panel of figure \ref{fig9}) remains weak during the first year but increases significantly afterwards. Like for natural photoperiod, the HF variations are not related to activity or inactivity phases.

\subsubsection{Comparison with experimental data}

The SPWVD analysis on simulated time-series with damped sinusoid $P_{spike}$ displays the same basic properties as the experimental data. In particular, there is a strong similarity between the LF components of the synthetic and experimental data. The fundamental pseudo-period of the damped sinusoid produces a circannual oscillation in the model outputs, that is detected as the LF activity by the SPWVD analysis.

The secretion rate $LH(t)$ generates high frequency spikes, that drive the fast oscillations in the solutions of (\ref{LHplasma}). For the chosen $P_{spike}(t)$ function and parameter values of (\ref{ParamSin}), the frequency window corresponding to this activity remains between 1 per 36 hours to 1 per 90 minutes. Thus, the HF component extracted by the SPWVD analysis results from this very high frequency activity distorted by the sampling process of twice per week frequency. Additionally, the HF components obtained from the experimental data and the model outputs display the same properties and seem to proceed from comparable underlying mechanisms.

\section{Discussion}
The present study gives convincing evidence for the existence of an endogenous circannual rhythm in Black Belly ewes. To our knowledge, this is the first time that such a rhythm is revealed in ewes originating from tropical regions. Up to now, it had only been studied in ewes originating from temperate regions. The SPWVD has been shown to be well suited both to extract this circannual hidden rhythm that was blurred by a very noisy environment, and to follow the changes all along three years. Since the SPWVD has not been already applied to neuroendocrine time series, we compared the parameters values given by SPWVD to another method used by Karsch {\it et al.} on five Suffolk ewes, which were exposed to the same experimental design including a LD-SD regimen (see \cite{Karsch89}). This protocol differed from ours on three points: (i) the Suffolk SD ewes originated from a temperate region (conversely tropical region for the Black Belly; however it was checked that estradiol plasma levels did not decrease throughout the experiment), (ii) they were studied during five years (conversely three years for the Black Belly), and (iii) the estradiol implant was replaced every year (conversely not replaced in Black Belly). In Karsch's study, the first phase was eliminated because it was induced by the initial transfer from long to short days. To compare our results to the ones of Karsch, we have averaged the parameter values over all the phases of LH secretion in each ewe, except the first one induced by the long day regimen. Despite these differences, the results of the two studies are in quite good agreement, so that they reinforce each other (table \ref{tab2}). The maximal LH level is similar between the two groups. The Black Belly ewes are characterized by a higher minimal LH amplitude, a longer LH phase duration and a shorter cycle duration. These results suggest that LH secretion is globally higher in the Black Belly ewes than in Suffolk SD. Nevertheless, despite these small differences, our results confirm those of Karsch {\it et al.}

\begin{table}
\caption{LH parameters values: comparison between Suffolk SD ewes (Karsch study) and Black Belly SD ewes (present study)}
\begin{center}
\begin{tabular}{|l|l|l|}
\hline
LH parameters& 
Suffolk SD& 
Black Belly\\
\hline
synchronisation (sem, day)& 
18& 
7 \\
\hline
phase duration (day)& 
126$\pm $14& 
132$\pm $4 \\
\hline
max. amplitude (ng/ml)& 
5.2& 
5.03$\pm $0.58 \\
\hline
min. amplitude (ng/ml)& 
0.4& 
1.33$\pm $0.30 \\
\hline
cycle duration (day)& 
336$\pm $ 16& 
285$\pm $15 \\
\hline
\end{tabular}
\label{tab2}
\end{center}
\end{table}

When applied to the model outputs obtained from the sinusoidal interspike interval function (\ref{PSin}), the SPWVD analysis correctly extracts the LF component, that is not impacted by the HF component. Additionally, the variations in the LF component allow to identify the circannual alternation from activity to inactivity phases. Hence the analysis of the model outputs confirms those obtained by the SPWVD analysis of the experimental data. More precisely, the synthetic and experimental data share the same basic features as far as the LF component is concerned: (i) the identification of four distinct phases, (ii) the increase in LF maximal amplitude, and (iii) the increase in LF basal level. The damped sinusoid seems to reproduce quite well the action of light on LH secretion in the experimental context.

Therefore, we can propose some assumptions about the response to photoperiod from the construction of function $P_{spike}$. It appears that, using a sinusoidal $P_{spike}$, we have obtained a typical LH level pattern comparable to what would be observed in a control ewe maintained under natural photoperiod. We are thus prone to interpret the response to natural photoperiod as sustained oscillations in the secretion frequency. In contrast, using a damped sinusoid $P_{spike}$ leads to the rhythmic priming fading progressively away, as a damped oscillator that is no more entrained. The sinusoidal part of the function represents the onset of the rhythmic process that may be due to the prior 3 month-exposure to an SD regime followed by a 3-month LD regime. During the first year cycle, the physiological response to the experimental photoperiod seems to be comparable to the response to a temperate photoperiod. Then, during the 3 year exposure to SD, only the endogenous rhythm can produce oscillations in LH frequency, and they fade away in the absence of suitable photoperiodic entrainment. This opens the way to study the photoperiod influence on the physiological system from the viewpoint of oscillator driving analysis.

Furthermore, we can infer the origin of the HF component of the experimental data from exploring which mechanisms produce the HF component of the synthetic data. Indeed, by construction, the solution of Equation 8 cannot display an oscillatory activity in the frequency window covered by the HF component. Nevertheless, with the damped sinusoid $P_{spike}$ function, the LH secretion function generates high frequency spikes (1 per 36 hours to 1 per 90 minutes), that drive the pulse frequency of the solutions of Equation 8. Hence, the HF component results from the distortion of the (very) high frequency pulsatility by the sampling process of twice per week frequency. The increase in HF amplitude ensues from this phenomenon: since the local maxima of $P_{spike}$ are getting smaller and smaller, the probability to obtain isolated spikes during expected inactivity phases becomes larger. We deduce that the experimental sampling on its own is sufficient to generate the HF components of increasing amplitude in the experimental LH time series. However, other factors may act as well: to distinguish between these artifactual features and an additional physiological mechanism, experimental sampling at a higher frequency should be undertaken.

\vspace{0.5cm}

\textbf{Acknowledgements:} The authors would like to thank Xochitl Hernandez for her technical assistance.

\vspace{0.5cm}

\textbf{Grants:} This work is part of the Large-Scale Project REGATE (REgulation of the GonAdoTropE axis)\footnote{http://www-rocq.inria.fr/who/Frederique.Clement/regate.html}.

\end{document}